\documentclass[a4paper,12pt, reqno]{amsart}
\usepackage[DIV=12, oneside]{typearea}
\usepackage[utf8]{inputenc}
\usepackage[T1]{fontenc}
\usepackage[english]{babel}
\usepackage{amsmath}
\usepackage{esint}
\usepackage{amssymb, amsthm, bbm, dsfont}
\usepackage{enumerate}
\usepackage{color, mathrsfs}
\usepackage{thmtools}
\usepackage{todonotes}
\usepackage{eso-pic}
\usepackage{isorot}
\usepackage[shortlabels]{enumitem}
\usepackage{microtype}
\usepackage[colorlinks=true, urlcolor=blue, linkcolor=blue,
citecolor=black]{hyperref}
\usepackage{array}
\usepackage{booktabs}
\usepackage{stmaryrd}

\theoremstyle{plain}

\declaretheoremstyle[bodyfont=\normalfont, qed=$\clubsuit$]{endsign}

\theoremstyle{definition}

\numberwithin{equation}{section}
\newcommand{\R}{  \mathbb{R}}

\newcommand{\cE}{  \mathcal{E}}
\newcommand{\cF}{  \mathcal{F}}

\renewcommand{\epsilon}{\varepsilon}

\renewcommand{\d}{\, \textnormal{d}}

\newlength{\widestlabel}
\setlist{itemindent=\parindent}
\setlist[enumerate]{labelwidth=\widestlabel, leftmargin=!}
\addto\extrasenglish{%

}

\author{Moritz Kassmann} 
\address{Fakultät für Mathematik\\
Universität Bielefeld \\
Postfach 100131 \\
D-33501 Bielefeld}
\email{moritz.kassmann@uni-bielefeld.de}

\title[On Dirichlet forms and semi-{D}irichlet forms]{On 
Dirichlet forms and semi-{D}irichlet forms \\ \vspace*{3ex}
{\scriptsize Remarks on the book \linebreak ``Semi-{D}irichlet 
forms and {M}arkov processes'' \linebreak by Yoichi Oshima}}

\subjclass[2010]{31-XX, 60-XX}
\keywords{Dirichlet forms, potential theory, Markov processes}
\thanks{The author would like to thank several participants of the Oberwolfach 
workshop ``Dirichlet Form Theory and its Applications '' in 2014 for helpful 
hints and discussions.}

\begin{document}

\begin{abstract}
One aim of this note is to give an 
overview of some developments in the 
area of Dirichlet forms. A second aim is to review the new book 
``Semi-{D}irichlet forms and {M}arkov processes'' by Yoichi Oshima. 
The book \cite{MR3060116}  appeared last year, but first versions were  
written as lecture notes 25 years ago. We first give a rather short and light 
introduction into the field of Dirichlet forms with a special 
emphasis on the subjects presented in the book under consideration. After 
a small  account on the history of Dirichlet forms we comment on the book by 
Oshima against the background of related works.
\end{abstract}

\maketitle

\section{Dirichlet forms and semi-Dirichlet forms} In short, the theory of 
Dirichlet forms is an advancement and an abstraction of potential theory. The 
theory of Dirichlet forms has created a lot of interesting and interrelated 
research since 1970. Dirichlet forms are related to probability theory, 
Riemannian geometry, pseudo-differential operators and mathematical physics. 
The strength and the beauty of this theory are that it provides a framework 
that connects spectral theory, functional analysis and stochastic processes in 
a natural way. Given some Hilbert space $L^2(X,m)$ of square-integrable 
functions on some topological space $X$ with measure $m$, a Dirichlet form is a 
pair $(\cE, \cF)$ of a bilinear form $(u,v) \mapsto \cE(u,v)$ for $u$ and $v$ 
from some domain $\cF \subset L^2(X,m)$. The domain $\cF$ itself, 
historically, is called Dirichlet space. Before discussing further requirements 
and examples, 
let us explain the main characteristics. A Dirichlet form is called symmetric 
if $\cE(u,v) = \cE(v,u)$ for all $u,v$. Whether $X$ is a locally compact 
separable metric space like a subset of $\R^d$ or whether $X$ is a more general 
infinite-dimensional state space has significant consequences for the whole 
theory. The book \cite{MR3060116} focusses on lower-bounded semi-Dirichlet forms 
which are defined on locally compact separable metric spaces $X$. 

Let us look at examples which are basic and were in the mind of those who 
founded the theory of Dirichlet forms. Assume $D \subset \R^d$ is open. We 
denote by $H^1(D)$ the space of all elements $u \in L^2(D)$ whose distributional 
derivatives $\nabla u$ are elements of $L^2(D)$. $H^1(D)$ is a Banach space with 
respect to the norm $\|u\|^2_{H^1(D)} = \|u\|^2_{L^2(D)} + \||\nabla 
u|\|^2_{L^2(D)}$. $H_0^1(D)$ denotes the closure of $C_c^\infty (D)$ with 
respect 
to $\|u\|_{H^1(D)}$. Set $\cE^{(1)}(u,v)=\int_D \nabla u \nabla v \d \lambda$ , 
$\cF^{(1)} = H^1(D)$, $\cF^{(2)} = H_0^1(D)$, where $\lambda$ denotes the 
Lebesgue measure. Then $(\cE^{(1)}, \cF^{(1)})$ and $(\cE^{(1)}, \cF^{(2)})$ are 
symmetric Dirichlet forms on $L^2(D, \lambda)$. For $\alpha \geq 0$ and $u,v \in 
\cF$ let us denote $\cE_\alpha(u,v) = \cE(u,v) + \alpha_0 (u,v)$ and 
$\cE_\alpha(u)=\cE_\alpha(u,u)$. The following conditions/properties turn out 
to be important. Assume $\alpha \geq 0$.
\begin{itemize}
 \item[($\cE$.1)] $\cE_\alpha(u) \geq 0$ for all $u \in \cF$.
 \item[($\cE$.2)] For some $K \geq 1$ and for all $u,v \in \cF$ 
$|\cE(u,v)| \leq K \sqrt{\cE_\alpha(u)} \sqrt{\cE_\alpha(v)}$. 
 \item[($\cE$.3)] For all $\beta > \alpha$ the domain 
$\cF$ is a Hilbert space with respect to the scalar product $\cE_\beta(u,v) + 
\cE_\beta(v,u)$.
 \item[($\cE$.4)] For all $u \in \cF$ the function $u^+ \wedge 1$ belongs to 
$\cF$ and $\cE(u^+ \wedge 1, u - u^+ \wedge 1) \geq 0$. 
\end{itemize}
These conditions define what is called a lower-bounded semi-Dirichlet form in 
\cite[Section 1.1]{MR3060116}. It makes sense that the author calls such forms 
Dirichlet forms in \cite{MR3060116} but it might be confusing for the reader of 
this review. The by now classical definition given in \cite{MR569058, MR2778606} 
requires $\cE$ to be symmetric, conditions ($\cE$.1), ($\cE$.3) to hold with 
$\alpha = 0$ and for all $u \in \cF$ the condition $\cE(u^+ \wedge 1, u^+ \wedge 
1) \leq \cE(u,u)$ which is stronger than ($\cE$.4). The notion of a nonsymmetric 
Dirichlet form is slightly more tricky. Again, one requires ($\cE$.1), ($\cE$.3) 
to hold with $\alpha = 0$. In this case one observes that ($\cE$.4) is 
equivalent to $\cE(u + u^+ \wedge 1, u - u^+ \wedge 1) \geq 0$. A nonsymmetric 
Dirichlet form now additionally requires 
\begin{itemize}
 \item[($\cE$.4')] For all $u \in \cF$ the function $u^+ \wedge 1$ belongs to 
$\cF$ and $\cE(u - u^+ \wedge 1, u + u^+ \wedge 1) \geq 0$. 
\end{itemize}
The above examples obviously are very special cases with additional features. 
First, they are symmetric forms. Second, we can choose $\alpha=0$ and $K=1$ by 
the Cauchy-Bunyakovsky-Schwarz inequality. ($\cE$.4) is called Markov property 
because it relates to the Markov property of the related stochastic process. 
Note that $\cE^{(i)}(u, u^+ \wedge 1) \geq \cE^{(i)}(u^+ \wedge 1)$ trivially 
for $i=1,2$. 

When relating these Dirichlet forms to extensions 
of $(-\Delta, C^\infty_c(D))$, readers might find 
\cite{Str96} quite informative despite the fact that the author does not hide 
his personal view and his preference for another approach. 

Let us provide an example which makes use of the flexibility of the above 
conditions. Define $\cE^{(3)}(u,v)=\cE^{(1)}(u,v) + \sum_{i=1}^d \int_D b_i 
\frac{\partial u}{\partial x_i} v \d \lambda$ for some functions $b_1, \ldots, 
b_d:D \to \R$ which either have bounded absolute values in $D$ or (for $d \geq 
3$) have the property that $\|\mathbf{b}\|_{L^d(D)}$ is finite and 
$\operatorname{div} \mathbf{b} $ is bounded from above. Here and below, we write 
$\mathbf{b}$ for the vector $(b_1, \ldots, b_d)^T$. The assumptions on 
$\mathbf{b}$ are tailored for conditions ($\cE$.1) and ($\cE$.2). This time, 
$\alpha \geq 0$ and $K \geq 1$ are chosen in dependence of $\mathbf{b}$. Thus 
the term {\it lower bounded} makes sense because of ($\cE$.1). The tuple 
$(\cE^{(3)}, H^1_0(D))$ is a lower bounded semi-Dirichlet form. In general, it 
is not a nonsymmetric Dirichlet form in the above sense. Note that our previous 
examples all were local forms in the sense that $\cE(u,v) =0$ if $u,v \in \cF$ 
have disjoint supports. There 
is a whole universe of nonlocal symmetric/nonsymmetric Dirichlet/semi-Dirichlet forms. 

There is a natural link between closed bilinear forms, semi-groups and 
resolvent operators. Assume $(\cE,\cF)$ satisfies ($\cE$.1)-($\cE$.3). Then 
there exist strongly continuous semigroups $(T_t)$, $(\widehat{T}_t)$ on 
$L^2(X,m)$ such that $\|T_t\| \leq e^{\alpha t}$, $\|\widehat{T}_t\| \leq 
e^{\alpha t}$, $(T_t f, g)=(f, \widehat{T}_t g)$ and for the resolvents 
$G_\beta, \widehat{G}_\beta$ given by $G_\beta f = \int_0^\infty e^{-\beta t} 
T_t f \d t$ and analogously for $\widehat{G}_\beta$: $\cE_\beta (G_\beta f, u) = 
(f,u) = \cE_\beta (u, \widehat{G}_\beta f)$ for all $f \in L^2(X,m)$, $u \in 
\cF$. The term {\it semi-Dirichlet form} relates to the fact that, different 
from $(T_t)$, the dual semi-group is not Markov in general.  

What we have explained so far, holds true in infinite dimensions, too. This 
changes, when it comes to the important concept of regularity. A lower bounded 
semi-Dirichlet form $(\cE, \cF)$ on $L^2(X,m)$ is called regular if $C_c(X)\cap 
\cF$ is (a) dense in $\cF$ with respect to the norm induced by $\cE_1(\cdot)$ 
and (b) dense in $C_c(X)$ with respect to the supremum norm. The concept of 
regularity needs to be changed significantly when working with infinite 
dimensional state spaces, which we comment on below. The major achievement of 
the theory of Dirichlet forms is that there is a correspondence 
between Hunt processes ($=$ quasi-left strong Markov processes) and regular 
Dirichlet forms if $X$ is a locally compact separable metric space. More 
precisely, there exists a Hunt process whose resolvent $R_\alpha f$ is a 
quasi-continuous modification $G_\alpha f$ for any $f \in L^\infty(X;m)$ and 
$\alpha>0$. Note that $R_\alpha f(x) = \int f(y) R_\alpha(x, \d y)$ where 
$R_\alpha(x, E) = \int_0^\infty e^{-\alpha t} p_t(x,E) \d t$ and $p_t$ is the 
transition function of the Hunt 
process. It is possible to give a complete characterization of all symmetric 
and nonsymmetric Dirichlet forms satisfying the sector condition in terms of 
right processes, see \cite{MR1214375}.

Given the relation between a given Dirichlet form and the corresponding Hunt 
process many properties of the form can be studied by investigating the process 
and vice versa. A fundamental result in the theory of regular symmetric 
Dirichlet forms is the formula of Beurling--Deny which provides a unique 
representation. Together with the results of Le Jan it leads to the following 
beautiful description which we give in a simple setting. Assume $X=D \subset 
\R^d$ is a domain and $(\cE, C_c^\infty(D))$ is a closable symmetric bilinear 
form on $L^2(D,m)$ satisfying the Markov property ($\cE$.4). Then $\cE$ can be 
expressed uniquely by
\begin{align*} 
\cE(u,v) &= \sum_{i,j=1}^d \int_D \frac{\partial u}{\partial_{x_i}} 
\frac{\partial v}{\partial_{x_j}} \nu_{ij}(\d x) \\
&\quad + \int_{D \times D} \big(u(y)\!-\!u(x)\big)\big(v(y)\!-\!v(x)\big) 
J(\!\d x \d y) + \int_D u(x) v(x) \kappa(\d x) \,, 
\end{align*}
where $\nu_{ij}, J$ and $k$ are nondegenerate positive Radon measures 
satisfying, among other properties, $\sum_{i,j=1}^d \xi_i \xi_j \nu_{ij}(K) \geq 
0$ and $\int_{K \times K} |x-y|^2 J(\!\d x \d y) < + \infty$ for any compact $K 
\subset \R^d$. From the probabilistic point of view, a major result in the 
theory of Dirichlet forms is the Fukushima decomposition of an additive 
functional into a martingale additive functional and an additive functional of 
zero energy. This decomposition is similar to the semi-martingale decomposition 
for Markov processes and leads to results which, in the simplest cases, can be 
obtained by the It\^{o} formula. We do not elaborate on this important result 
here. 

The development of Dirichlet forms and corresponding Markov processes on 
infinite-dimensional state spaces is motivated by questions arising in quantum 
field theory and interacting particle systems. The main mathematical challenge 
is to find a substitute for the notion of regular Dirichlet forms. To this end, 
the notion of quasi-regularity for Dirichlet forms is introduced. Leaving 
technicalities aside, a Dirichlet form is quasi-regular if and only if it is 
quasi-homeomorphic to a regular Dirichlet form on a locally compact metric 
space. It can be shown that a Dirichlet form is associated with a nice Markov 
process if and only if it is quasi-regular. The quasi-homeomorphism allows the 
results developed for regular Dirichlet forms to be applied to quasi-regular 
Dirichlet forms on general infinite-dimensional state spaces. Note that the term 
``quasi'' relates to exceptional sets which are defined with respect to 
capacity. These exceptional sets appear naturally and are abundant in the theory 
of Dirichlet forms.
 One task is to overcome them using regularity theory.

\section{A small and incomplete account on the history of Dirichlet forms} The 
articles \cite{MR0106365, MR0098924, MR0284609} study the domains of Dirichlet 
forms as Hilbert spaces and show that the concept of Dirichlet space captures a 
lot of the classical potential theory. Usually, these studies are regarded as 
the birth of the analytic side of Dirichlet form theory. The correspondence with 
a strong Markov process is provided in \cite{MR0295435}. This review is not the 
right place to list all articles that contributed to this theory. We restrict 
ourselves to monographs and to those articles which are related to the focus of 
\cite{MR3060116}. 

The books \cite{MR0386032, Fuk75, MR0451422, MR569058} lay out the foundations 
of symmetric Dirichlet forms and corresponding Markov processes. It is important 
to note that nonsymmetric forms and related stochastic calculus were studied 
already at the very beginning of the theory, e.g. in \cite{MR0270445, MR0399491, 
MR0399490, MR0457756, MR0588389, MR0386030,  MR571671, MR0467934, MR0467934, 
MR508949, MR551635, MR506612, MR661618, MR909022}. Note that these works benefit 
from the corresponding theory for elliptic differential operators in divergence 
form worked out in \cite{MR0192177}. A standard reference for nonsymmetric forms 
is the monograph \cite{MR1214375} which appeared roughly at the same time as 
\cite{MR1133391} and the first edition of \cite{MR2778606} which is an extension 
of \cite{MR569058} and nowadays is the main reference for the field. Note that 
the main emphasis of \cite{MR1214375} is on the development of the theory of 
Dirichlet forms on general state spaces. As 
mentioned above, the 
motivation to relax the assumptions (locally compactness) on the state space is 
closely connected to mathematical physics. First studies in this direction 
include \cite{MR0455133, MR0461153}. Since the book under review does not add 
results in this direction we do not provide more references on this important 
development and refer the interested reader to the discussion in 
\cite{MR1038449, MR1133391, MR1214375, MR1200639, MR1632609, 
MR1772835, MR2009816, AMR14}. 
Note that the exposition in \cite{MR1632609} goes beyond the scope of other 
books and covers truely nonsymmetric (without sector condition) and rather 
general time-dependent Dirichlet forms for infinite-dimensional state spaces.   

A second field of current interest which is not touched by the book under 
consideration is the connection between geometry and Dirichlet forms. Again, we 
decide not to list many articles but rather give only a few hints where to find 
more information. The proceedings volume \cite{MR1652277} might be a good start 
because it contains several related articles. On the one hand, Dirichlet forms 
provide a tool to define a Laplacian on general state spaces. On the other 
hand, they provide a framework for results which are robust with 
respect to geometric quantities.  The geometric significance 
of Harnack inequalities and heat kernel bounds for Dirichlet forms are studied 
in \cite{MR1301456, MR1378473, MR1355744, MR1387522}, see also the references 
therein. In typical situations, the Gaussian short time off-diagonal 
behavior of the heat kernel is a function of the intrinsic 
distance. This holds true for general strongly local Dirichlet forms 
\cite{MR1988472, MR2176031}. Aronson-type bounds have 
been studied in metric measure spaces and on fractals using the theory of 
Dirichlet forms. The works \cite{MR1668115, MR1840042, MR2228569, MR2743439, 
MR2639315, MR2998918, MR2919892, MR2962091, MR2970465, MR3250370, MR3194164, 
AnBa14} contain several important results. The recent book 
\cite{MR3155209} is a good source for the relation of functional inequalities 
and local Dirichlet forms 
in a general context. Note that the aforementioned contributions mainly 
concentrate on local Dirichlet forms. Similar studies for nonlocal Dirichlet 
forms and their relation 
to geometry seem to be much more subtle. Let us also mention that Dirichlet 
forms can be applied to other areas like discrete goups or random media. Two 
chapters of \cite{MR1218884} address discrete groups and estimates of the decay 
of convolution powers of probability measures on these groups. For applications 
to random media see \cite{Kum14}.

Last, let us mention the monograph \cite{MR2849840}. It provides a 100-pages 
summary of the theory of symmetric regular and quasi-regular Dirichlet forms 
which can be used as a first read. Moreover, it treats new developments in more 
specialized fields, i.e., trace processes, boundary theory and reflected 
Dirichlet spaces for regular symmetric Dirichlet forms. 

As explained above, nonsymmetric Dirichlet forms were studied right on from the 
beginning of research activities around Dirichlet forms. The lectures of Y. 
Oshima at Friedrich-Alexander-Universit\"at Erlangen-N\"urnberg in 1988 and 1994 
contributed significantly to this development. It is not clear to the author of 
this review when the notion of a semi-Dirichlet form was used for the first 
time. Regularity resp. quasi-regularity of semi-Dirichlet 
forms are studied in \cite{MR1323103, MR1341121}. Several examples of 
semi-Dirichlet forms can be found in \cite{MR1319695, MR2952094, MR3273873}. 
Chapter 7 
of \cite{MR2153655} contains several results on quasi-regular semi-Dirichlet 
forms. Recently, the Fukushima decomposition and the Beurling-Deny formulae have 
been studied for semi-Dirichlet forms \cite{MR2057023, MR2990583, 
MR2994113, MR2584980, MR2920204}. See also the very recent survey 
\cite{MSW14} and which the author is  thankful for having been 
provided with.

\section{Remarks on the new book by Yoichi Oshima} The new book \cite{MR3060116} 
by Y. Oshima extends the theory of Dirichlet forms on locally compact separable 
metric spaces. After having set up the analytic theory of lower bounded 
semi-Dirichlet forms in the first two chapters, the relation to Hunt processes 
is studied in Chapter 3. Chapters 4 and 5 are devoted to additive functionals 
and decompositions. Thus the book can be viewed as a natural extension of 
\cite{MR2778606}. The book is carefully written and has got a nice layout. It 
certainly will serve as a standard reference for its respective 
field. Because of the development of the theory of 
Dirichlet forms during the past 25 years, it is natural that there are by now 
other monographs. The standard reference for symmetric Dirichlet forms on 
locally compact separable metric spaces remains \cite{MR2778606}. Symmetric 
Dirichlet forms in a more general framework are presented nicely also in 
\cite{MR2849840}. If interested  in symmetric or nonsymmetric Dirichlet forms 
satisfying the sector condition on more general state spaces, then 
\cite{MR1214375} is the first choice. A special feature of \cite{MR3060116} is 
the treatment of time-dependent Dirichlet forms in Chapter 6 which is an 
advancement of parabolic potential theory. These questions are also covered in 
\cite{MR1632609} and even in much greater generality. However, the presentation 
in \cite{MR3060116} is easier accessible and sufficient for many purposes.  
Altogether, the new book \cite{MR3060116} is a welcome addition 
to the literature. The accuracy of the presentation and the similarity of the 
style to the one of \cite{MR2778606} will be appealing to many reasearchers. As 
mentioned above, the main source for the book are the lectures by the author on 
nonsymmetric Dirichlet forms from 1988 and 1994 in Erlangen. Although the 
corresponding lecture notes are not available to the public, they have been 
circulated with the author's permission among interested colleagues. In this 
way, the book certainly has had an impact on the theory long before its 
publication.

{\small 

\newcommand{\etalchar}[1]{$^{#1}$}

}
\end{document}